\documentclass[11pt]{amsart}

\usepackage{amsmath, amssymb, amsbsy, amsfonts, amscd, epic, eepic}

\numberwithin{equation}{section}

\newtheorem{theorem}{Theorem}[section]  
\newtheorem{definition}[theorem]{Definition}  
\newtheorem{example}[theorem]{Example}  

\newtheorem{lemma}[theorem]{Lemma}  
\newtheorem{proposition}[theorem]{Proposition}  
\newtheorem{corollary}[theorem]{Corollary}  
\newtheorem{remark}[theorem]{Remark}  
\newtheorem{conjecture}[theorem]{Conjecture}

\def\cal{\mathcal}  

\newcommand{\A}{{\cal A}} 
\newcommand{\B}{{\cal B}} 
\newcommand{\C}{{\cal C}}

\newcommand{\Lo}{{\cal L}}

\newcommand{\N}{{\cal N}}

\newcommand{\V}{{\cal V}}

\newcommand{\X}{{\cal X}}

\def\\N{\mathbb{N}}
\def\CC{\mathbb{C}}

\def\PP{\mathbb{P}}
\def\ZZ{\mathbb{Z}}

\def\P{\mathbb{P}} 

\font\tengoth=eufm10  \font\fivegoth=eufm5
\font\sevengoth=eufm7
\newfam\gothfam  \scriptscriptfont\gothfam=\fivegoth
\textfont\gothfam=\tengoth \scriptfont\gothfam=\sevengoth

\title{Realization of finite Abelian groups by nets in $\PP^2$}
\author{Sergey Yuzvinsky}
\email{yuz@math.uoregon.edu}
\address{University of Oregon\\
Eugene, OR 97403}
\subjclass{Primary: 52C30, 05B30; Secondary: 14H50}
\date{\today}
\begin{document}
\begin{abstract}
In the paper, we study special configurations of lines and points
in the complex projective
plane, so called $k$-nets. We describe the role of these configurations in studies
of cohomology on arrangement complements.
Our most general result is the restriction on $k$ -
it can be only 3,4, or 5. The most interesting class of nets is formed by 3-nets
that relate to finite geometries, latin squares, loops, etc. All known examples
of 3-nets in $\PP^2$ realize finite Abelian groups. We study the problem what
groups can be so realized. Our main result is that, except for groups with all
invariant factors under 10,
realizable groups are isomorphic to subgroups of a 2-torus. This follows from
the `algebraization' result asserting that in the dual
plane, the points dual to lines of a net lie on a plane cubic.
\end{abstract}

\maketitle
\bigskip
\section{Introduction}
The notion of a net, that is a discrete analogue of a web, has a long history
related to discrete geometries, latin squares, loops, etc. (e.g., see \cite
{DK}). It turned out
recently that nets play a special role in the study of the cohomology of local
systems on the complements of complex line arrangements (see section 2).  Since the
latter discovery attracted the author to nets, we do not even define
abstract nets in the paper but from the beginning restrict our considerations
to nets in
the complex projective plane $\PP^2$. These nets can also be viewed
as complex representations of certain matroids of rank 3.

\begin{definition}
\label{def_nets}
Let $k$ be an integer, $k\geq 3$. A {\it $k$-net} in $\P^2$ is a pair $(\A,\X)$
where $\A$ is a finite set 
of lines partitioned into $k$ subsets $\A=\bigcup_{i=1}^k\A_i$ and $\X$ 
is a finite set
of points subject to the following conditions:

(i) for every $i\not= j$ and every $l\in\A_i, \ l'\in\A_j$, we have $l\cap
l'\in \X$;

(ii) for every $X\in\X$ and every $i$ ($i=1,2,\ldots,k$) there exists a unique
$l\in\A_i$ passing through $X$.
\end{definition}
The definition of a $k$-net implies
immediately that $|\A_i|=|l\cap\X|$ for every
$i, \ 1\leq i\leq k,$ and every $l\in\A$.
This number, which we denote by $m$, is called {\it the order of the net}.
We will often use the term $(k,m)$-net for a $k$-net of order $m$.
For a $(k,m)$-net we have $|\X|=m^2$.
We will always consider nontrivial nets defined by the condition $m>1$.

Notice that we do not put any restrictions on the combinatorics of lines inside a
class; they can form an arbitrary arrangement from general position one to
a pencil.

Our first result about nets (Theorem \ref{euler}) asserts that the only
possible values for $k$ are 3,4, and 5 (cf. \cite{LY}). 
By putting extra conditions on combinatorics inside classes one can
restrict $k$ even more.
For instance, if two classes are pencils then $k=3$.

Since 3-nets form the most interesting class of nets and only one example of
other nets in $\PP^2$ is known to me, we restrict our considerations to $k=3$.
In this case, there is the pairing $\A_1\times\A_2\to\A_3$ define by assigning
to $l_i\in\A_i$ ($i=1,2$) the line from $\A_3$ passing through $l_1\cap l_2$.
After any identification of each class with a set $H$ this paring defines
a binary operation on $H$ converting $H$ into a quasigroup. A general
theorem of quasigroup theory (e.g., see \cite{DK}, Theorem 8.1.4) asserts that
one can always find identifications such that $H$ is a loop, that is all
group axioms hold for it but perhaps the associativity. If this loop is a group
we say that the net {\it realizes} the group. 

There are conditions on a net in $\PP^2$ implying that it represents a
particular group. For instance if all classes of a $(3,m)$-net
are pencils then the net realizes
$\ZZ_m$ (Proposition \ref{prop_pencils}). 
No example of a 3-net in $\PP^2$ is known to me where the net does not
realize an Abelian group. The main part of the paper deals with the question
which Abelian groups are realizable by nets in $\PP^2$.

Our main result in this direction is somewhat similar to the theorem of 1924
by H.Graf and R.Sauer about algebraization of webs. Following theory of webs we
say that a 3-net is {\it algebraic} if there exists a cubic curve in the dual
plane $\PP^*$ passing through all the points dual to the lines of the net.
We prove (Theorem \ref{thm_cyclic} and Theorem \ref{thm_oncubic}) that
a nets of a wide class are algebraic. This class contains the
nets realizing arbitrary cyclic groups and also Abelian groups 
with elements of order larger than 9.

Using the algebraization result we can prove our negative
results about realization. The main statement (Corollary
\ref{nonrealization}) asserts that a finite Abelian group $H$ 
with an element of
order larger than 9 can be realized if and only if 
it is the direct sum of not more than two cyclic groups. The `if' part is true
without any restriction on the orders of elements (Theorem \ref{cubic}).
If $H$ has an element
of order larger than 6 then it can be realized by a net with one class being a
pencil if and only if $H$ is cyclic. Also any cyclic group can be realized with
all three classes being pencils and this realization is projectively unique
(Proposition \ref{prop_pencils}).

We conjecture that the realization result holds without any condition on the
order of the elements. Also we suggest problems to find more examples of 
4- and 5-nets in $\PP^2$ and 3-nets not realizing Abelian groups (or to prove
that such examples do not exist).

The paper is organized as follows. Section 2 is devoted completely to
motivations and is not needed for understanding the rest of the paper.
In section 3, we prove general restrictions on $k$ for a $k$-net in $\PP^2$. In
particular the complete classification of nets whose 3 classes are pencils is
given. In section 4 (and for the rest of the paper), we restrict our
considerations to 3-nets and study which Abelian groups can be realized by them. 
We exhibit first positive and negative examples of realization. In
section 5, we prove our main theorems about the property of nets being
algebraic
and about Abelian group realization. For that we need first to study the
operations defined by the collinearity of points on reducible cubics. Finally in
section 6, we state conjectures and unsolved problems.

\section{Motivation: essential components of characteristic varieties}
\bigskip
As in Introduction, we denote by $\A$ a finite set $\{l_1,\ldots,l_n\}$ of
lines in the complex projective plane $\PP^2$ and put $M=\PP^2\setminus \bigcup
_{i=1}^nl_i$. 
The algebra $H^*(M,\CC)$ is naturally isomorphic
to the Orlik-Solomon algebra $A$ of $\A$. The latter algebra is generated by
skew-commutative degree 1 generators $e_1,\dots, e_n$ satisfying certain extra
relations of degree at least 2 (see \cite{OT}). 
Computing the cohomology of local systems on $M$ is a much more subtle task.
Every rank one local system on $M$ is defined by a character
$H_1(M)\to \CC^*$ whence the set of all those systems form the $n$-dimensional
torus $T=H^1(M,\CC^*)$. 
For every $k\geq 1$ we put
$\Sigma_k=\{\Lo\in H^1(M,\CC^*)|\dim H^1(M,\Lo)\geq k\}$. This is an algebraic
subvariety of the torus called {\it the $k$th characteristic variety of $M$}
(or of $\A$).

Here we will consider only the irreducible components 
of $\Sigma_k$ containing the identity {\bf 1} =$(1,1,\ldots,1)$ of the torus,
i.e., the constant local system. For such a component it is possible to describe
its tangent space at {\bf 1}. 
Note that the multiplication by an
element $a\in A_1$ defines a differential on $A$ of degree 1 whence the
cohomology spaces $H^p(A,a)$. This allows one to define so called {\it resonance
subvarieties} of $A_1$ via $R_k=\{a\in A_1|\dim H^1(A,a)\geq k\}$. 
Now since the above generators of $A$ form a basis of $A_1$ we can identify
$A_1$ with the tangent space to $T$ at {\bf 1} (i.e., the Lie algebra of $T$). 
 The following result is just a little stronger than the results proved
in \cite{LY}.

\begin{theorem}
\label{components}
 For every component $\V$ of $\Sigma_k$ containing {\bf 1}
its tangent space is a component of $R_k$. 
\end{theorem}

Now let us summarize results from \cite{LY} (discussed with more detail
in \cite{Yu}, Ch. 7). Every irreducible component of $R_k$ is defined by
a subarrangement $\B\subset\A$ and a set $\X$ of multiple points of $\B$.
If $J$ is the incidence matrix of the pair $(\X,\B)$ and $E$ is the $|\B|\times
|\B|$-matrix whose all the entries equal 1 then the matrix $Q=J^{t}J-E$
should have at least $k+2$ indecomposable components of affine type (see
\cite{Ka}) that cover every $X\in\X$
and $\B$ should be maximal with respect to this property.
In the rest of the paper we will consider only {\it essential} components for
which $\B=\A$.

One of the corollaries of the description of a component of $R_k$ is that
the multiplicity of every point from $\X$ is at least $k+2$ (\cite{LY},
Proposition 2.5). Thus the most `economical' components are those for which
all points of $\X$ have the same multiplicity equal to $k+2$. These components
correspond to the full graphs in the sense of \cite{LY}, Sect. 6, and the
combinatorics of them has been described in \cite{Yu}. 
In terms of nets in $\P^2$ (Definition \ref{def_nets}) the latter
result just asserts
that those components correspond to the case where
the pair $(\A,\X)$ is a $k+2$-net.

Given a $k+2$-net $(\A,\X)$ in $\P^2$ it is easy to write explicitly
the respective essential component $V=V(\X)$ of $R_k$ for $\A$. The partition 
$\A=\bigcup_{i=1}^{k+2}\A_i$ induces the respective partition of the basis
$(e_1,e_2,\ldots,e_n)$ of $A_1$. Then $V$ consists of all $a=(x_1,\ldots,x_n)\in
A_1$ such that $\sum_{i=1}^nx_i=0$ and $x_i=x_j$ on each block of the partition.

A subset $\B\subset \A$ formed by a pencil of lines also defines a component of
$R_k$, $k=|\B|-2$. This component is not essential unless $\B=\A$ and is called
local. The local component corresponds to the trivial $(k,1)$-net formed by $\B$
and the base of the pencil. 
\begin{example}
\label{braid}

The arrangement of lines in $\P^2$ given by the linear functionals
$(x,y,z,x-y,x-z,y-z)$ can be partitioned into 3 parts $\{x,y-z\},\{y,x-z\},$
and $\{z,x-y\}$. It is easy to check that this partition defines a (3,2)-net
with all four triple points of the intersections. The respective component
of $R_1$ is the unique non-local (essential) component.

\end{example}

\section{Restrictions on $k$-nets in $\P^2$}

Unlike theory of the abstract $(k,m)$-nets 
(see \cite{DK}) where there are no general
restrictions on $k$ and $m$, one can obtain quite a few restrictions on these
parameters for nets in $\P^2$.

First we observe that a $(k,m)$- net induces a pencil of 
curves of degree $m$ with $k$ very
special fibers. (In \cite{LY} these pencils where used in more general
situation using deep results by Arapura; here we introduce them elementarily.)

For a $(k,m)$-net $(\A,\X)$ with a partition $\A=\bigcup_{i=1}^k\A_i$
consider $k$ split curves of degree $m$ given by
$C_i=\bigcup_{l\in\A_i}l$, $1\leq i\leq k$. The defining ideal of $C_i$ is
generated by the polynomial $g_i=\prod_j\alpha_{ij}$ where 
$\{\alpha_{ij}|j=1,\ldots,m\}$ are some
fixed linear functionals whose kernels runs through the lines from $A_i$.

The following lemma is a simple case of Max Noether's $AF+BG$ Theorem.
\begin{lemma}
\label{pencils}
All curves $C_i$ are fibers of the pencil generated by any two of them.
\end{lemma}
\proof
Choose two of the curves, say $C_1$ and $C_2$,
Since the curves intersect
at $m^2$ (simple) points 
the ideal $I$ generated by $g_1$ and $g_2$ in the polynomial
ring is radical. Also any other curve $C_i$ passes through all these intersection
points. By the Nullstellensatz we have $g_i\in I$. 
Since the degree of $g_i$ is again $m$ we obtain the result.
\qed

\medskip
Lemma \ref{pencils} implies that there is a rational map 
$\phi:\P^2\to \P^1$ with poles at
the points of the net that sends $p\in\P^2$ to the point $(c:d)\in\P^1$ such that
$cg_1(p)+dg_2(p)=0$. Blowing up all the points of the net we obtain a variety
$P'$ and a regular map $\psi:P'\to \P^1$ whose preimages of points are the
proper transforms of fibers of the pencil. This allows us to prove the following 
properties of the net (cf. \cite{LY}, Proposition 7.3).

\begin{theorem}
\label{euler}
(i) For an arbitrary $(k,m)$-net in $\P^2$ the only possible values for $(k,m)$
are: $(k=3, m\geq 2),(k=4, m\geq 3), (k=5, m\geq 6)$.

(ii) If one class of the net is a pencil then the only possible value are:
$(k=3, m\geq 2)$ and $(k=4, m\geq 4)$.

(iii) If two classes are pencils then the only possible value
for $k$ is 3.
\end{theorem}
\proof
To cover all the cases let us assume that $r$ classes of the net are pencils,
$0\leq r\leq k$. 
We estimate the Euler characteristic $\chi(P')$. Since we got $P'$ by blowing up 
$m^2$ points of $\P^2$ we have $\chi(P')=3+m^2$. On the other hand, let us
estimate the Euler characteristics $\chi$ of fibers of $\psi$. The general
position fiber is  a smooth curve of degree m whence $\chi=3m-m^2$. At least
$k$ fibers are unions of $m$ lines whence for each of these fibers 
$\frac{m(5-m)}{2}\leq \chi\leq m+1$ where the left hand side value is reached 
by generic configurations of lines and the right hand side value 
by a pencil of lines. Combining
this information we have
$$3+m^2\geq (2-k)(3m-m^2)+ r(m+1)+(k-r)\frac{m(5-m)}{2}\eqno(1.1)$$
or equivalently (since $m>1$)
$$k\leq 6\frac{m-1}{m}-r\frac{m-2}{m}.\eqno(1.2)$$
Ignoring the second summand we immediately get the first statement.

If $r=1$ the inequality (1.2) becomes 
$$k\leq 5-\frac{4}{m}$$
which implies the second statement.
Finally if $r=2$ we have 
$$k\leq 4-\frac{2}{m}$$
which implies the last statement.                     \qed

\medskip
It turns out that nets whose all classes are pencils can be 
classified completely up to projective isomorphism.
\begin{proposition}
\label{prop_pencils}
For any $(3,m)$-net $N=(\A,\X)$ in $\P^2$ 
whose all classes are pencils there exists a homogeneous coordinate system
such that the classes $\A_i$ of $N$ are given by the equations
$x^m-y^m=0, \ x^m-z^m=0$ and $y^m-z^m=0$ respectively.
In particular $N$ realizes $\ZZ_m$.
\end{proposition}
\proof
Since each class $\A_i$ is a pencil of lines it defines a point $P_i$,
the base of the pencil. Consider two possible cases.

{\it Case 1}. The points $P_i$ are collinear, lying on a line $\ell$. 
Deleting $\ell$ from $\P^2$ we
reduce this case to a 2-net in $\CC^2$ whose each class consists of lines
parallel to each other (cf. \cite{DK}, 8.3). Choose an affine coordinate system
so that the first class is formed by the lines $x=a_i$, the second one by $y=b_i$,
and the third one by $x+y=c_i$, where $a_i,b_i,c_i\in\CC$ and $i=1,2,\ldots, m$.
We can choose the origin of the system so that the lines 
$x=0$, $y=0$, and $x+y=0$ are in the respective classes.
Now we identify each class with a subset of
$\CC$ assigning to each line the right hand side of its equation. 
Under this identification the resulting quasigroup (see Introduction)
becomes a subgroup of
the additive group $\CC$ which contradicts the finiteness of $N$.

{\it Case 2}. The points are not collinear.
Choose a coordinate system such that $P_1=(0:0:1),
P_2=(0:1:0), P_3=(1:0:0)$. This allows us to write equations of lines in the
classes as $a_ix-y=0$, $x-b_iz=0$, and $y-c_iz=0$ respectively where
$a_i,b_i.c_i\in\CC^*$ and $i=1,2,\ldots,m$. Again we can normalize the equations
of one line from each class to $x-y=0$, $x-z=0$, and $y-z=0$
respectively and then identify classes with $\CC^*$ assigning to a line the
coefficient $a_i$, $b_i$, or $c_i$ respectively.
We obtain then a quasigroup which is a subgroup of $\CC^*$ whence isomorphic to
$\ZZ_m$. Fixing a primitive $m$th root of 1 as a generator and multiplying
equations in each class we obtain the result.                   \qed

\begin{remark}\label{rmk_pencils}
 I have a little doubt that Proposition \ref{prop_pencils} is
known. Section 8.3 of \cite{DK} is devoted to similar topic. However I could not
find there the statement we need. The proof is easy and is included
for convenience of the reader.
\end{remark}

In the next section we will need the following corollary of  Proposition
\ref{prop_pencils}.
\begin{corollary}
\label{cor_pencils}
Let $N=(\A,\X)$ be a $(3,mn)$-net $(m\geq 3)$
 such that each class $\A_i$ is partitioned into
$n$ blocks of size $m$, i.e., $\A_i=\cup_{j=1}^n\A_{ij}$, and
for every two $i,j$ there is a $k$ such that $\A_{1i},\A_{2j}$, and $\A_{3k}$
form 3 classes of a $(3,m)$-subnet $N_{ij}$ of $N$. Then if for every 
$i,j$ each class of $N_{ij}$ is a pencil, every class of $N$ is also a pencil.
\end{corollary}
\proof
Denote by $X_{ij}$ the base of the pencil $\A_{ij}$ and 
suppose the conclusion is false. Then without loss of generality we can assume
that $X_{11}\not=X_{12}$. Consider the two nets $N_{11}$ and $N_{21}$.
They have the same second class $\A_{21}$ and distinct first classes $\A_{11}$ and
$\A_{12}$. By Proposition \ref{prop_pencils} there is a projective isomorphism
$\phi$ of $\P^2$ fixing every line from $\A_{21}$ and mapping lines from
$\A_{11}$ to lines from $\A_{12}$. Since $m\geq 3$, $\phi$ fixes every line
passing through $X_{21}$, in particular the lines $l_1=(X_{21}X_{11})$ and
$l_2=(X_{21}X_{12})$. (Here and below
$(ab)$ denotes the line passing through points $a$
and $b$). Since besides $\phi(X_{11})=\phi(X_{12})$ we see that
the points $X_{11},X_{12},X_{21}$ are collinear. 

Now there exists $k$ such that
$\A_{3k}$ is the third class in $N_{11}$. Substitute the class
$\A_{3k}$ for $\A_{21}$ and consider $N_{11}$ together with the $(3,m)$-subnet
defined by $\A_{12}$ and $\A_{3k}$. We see similarly to the above
that $X_{11},X_{12},X_{3k}$
are collinear also. This implies that the bases $X_{11},X_{21},X_{3k}$ of the
pencils of the net $N_{11}$ are collinear which contradicts Proposition
\ref{prop_pencils}. The contradiction completes the proof.           \qed

\medskip
\begin{example}
\label{hessian}
The only known to me example of a $4$-net in $\PP^2$ is the Hessian
configuration. The set of points for the net consists of the nine inflection
points of a non-singular plane cubic and the set of lines consists of the 12
lines each passing through three of these points. In this example $m=3$
and each of the four classes is in general position (i.e., its lines intersect
at three distinct points).
\end{example}

\medskip
No example of a $5$-net in $\PP^2$ is known to me.
\bigskip
\section{Realization of finite Abelian groups; first examples}
\bigskip
In the rest of the paper we will focus our attention on 3-nets.  
Let us recall the important feature of 3-nets from Introduction.
For a 3-net any pair $(l_1,l_2)$, $l_i\in\A_i$ ($i=1,2$) defines a unique
line $l_3\in\A_3$ passing through the point $l_1\cap l_2$. Thus we obtain a
pairing $\A_1\times \A_2\to\A_3$ or, after some identification of each of three
classes with a set $H$, a binary operation on $H$. This operation defines on $H$ a
structure of a quasigroup or (under appropriate identification) even a loop
(see \cite{DK}). If there exists an identification such that this loop is an
Abelian group then we say that this group is realized by the 3-net.

We remark at this point that all known examples of 3-nets in $\P^2$ realize
Abelian groups. So the rest of the paper is devoted to the question which
finite Abelian groups can be realized by a 3-net (in $\P^2$).

First we prove a negative result.

\begin{lemma}
\label{Z_2 pattern}
Suppose there are 9 points $a_1,a_2,a_3,b_1,b_2,b_3,c_1,c_2,c_3$ such that
the following triples are collinear:
$\{a_1,b_1,c_1\},\{a_1,b_2,c_2\},\{a_1,b_3,c_3\},$\hfill
\break $\{a_2,b_1,c_3\},\{a_2,b_3,c_1\},
\{a_3,b_2,c_3\}$ and $\{a_3,b_3,c_2\}$. Then the lines \hfill\break
$\{(a_2a_3),(b_1b_2),(c_1c_2)\}$ intersect at one point.
\end{lemma}
\proof The triangles $b_1b_2c_3$ and $c_1c_2b_3$ form a Desargus configuration,
more precisely the lines passing through pairs of corresponding vertices
intersect at $a_1$. Besides $(b_1c_3)\cup (c_1b_3)=\{a_2\}$ and 
$(b_2c_3)\cap (c_2b_3)=\{a_3\}$. Now the result follows from 
the second Desargus theorem.                \qed

\begin{theorem}
\label{Z_2}
The group $\ZZ_2^3$ cannot be realized.
\end{theorem}
\proof
Suppose that $\ZZ_2^3$ is realized
by 3 sets of points $\{a_i\},\{b_i\}$ and $\{c_i\}$ where $i=1,2,\ldots,8$.
Let us consider the matrix $C$ of the pairing
generated by the collinearity relation.  For that we enumerate the rows of this
matrix by $a_i$ and the columns by $b_j$. The ($a_i,b_j$)-entry is the point
$c_k$ lying on the line $(a_ib_j)$. Without loss of generality we can assume
that
$$C=\left(\begin{matrix}A&B\\
                        B&A\end{matrix}\right)$$
where 
$$A=\left(\begin{matrix} c_1&c_2&c_3&c_4\\
                         c_2&c_1&c_4&c_3\\
                         c_3&c_4&c_1&c_2\\
                         c_4&c_3&c_2&c_1\end{matrix}\right)$$
and
$$B=\left(\begin{matrix} c_5&c_6&c_7&c_8\\
                         c_6&c_5&c_8&c_7\\
                         c_7&c_8&c_5&c_6\\
                         c_8&c_7&c_6&c_5\end{matrix}\right).$$
Using Lemma \ref{Z_2 pattern} for rows $a_1,a_3,a_4$ and columns $b_1,
b_2,b_4$ we learn that
the lines $\{(a_3a_4), (b_1b_2),(c_1c_2)\}$ intersect at a common point, say $x$.
Using the lemma for rows $a_1,a_5,a_6$ and columns $b_1,b_2,b_5$ we learn
that $x\in (a_5a_6)$ and similarly $x\in (a_7a_8)$. Using the lemma again for
rows $a_5,a_7,a_8$ and columns $b_5,b_6,b_7$ we have $x\in (b_5b_6)\cap (a_1a_2)$. 

Now consider the quadrangle $a_1b_1b_2a_2$. We have $(a_1b_1)\cap
(a_2b_2)=\{c_1\}$ and $(a_1b_2)\cap (a_2b_1)=\{c_2\}$. Thus by the Diagonal theorem
the points $x_1=(c_1c_2)\cap (a_1a_2)$ and $x_2=(c_1c_2)\cap (b_1b_2)$ 
form a harmonic
pair for the pair $\{c_1,c_2\}$. By the inclusions in the previous paragraph
however $x_1=x=x_2$. The contradiction completes the proof.     \qed

\medskip
Using similar (but simplier) argument one can proof the following theorem.

\begin{theorem}
\label{Z_2 pencil}
There is no representation of $\ZZ_2^2$ with at least one class being a pencil.
\end{theorem}

\medskip
Now we prove a positive result about realization of finite Abelian groups
by 3-nets.

\begin{theorem}
\label{cubic}
Let $H$ be a finite subgroup of a two dimensional torus.
Then there exists a 3-net in $\P^2$ realizing $H$.
\end{theorem}
\proof
Fix a non-singular plane cubic curve $C$ and convert it to the two dimensional
torus choosing one of the flexes of $C$ as the neutral element. Then
three points $p_1,p_2,$ and $p_3$ from $C$ are collinear if and only if
$p_1+p_2+p_3=0$. We can identify $H$ with a subset of $C$ and put $m=|H|$. 
Now choose
$\alpha,\beta\in C/H$ such that $\alpha,\beta$, and $-\alpha-\beta$ are all
distinct. 
The union of these three cosets is a set of $3m$ points 
partitioned in three classes. It is clear that the lines in the
projectively dual plane dual to these points 
form a 3-net and this net realizes $H$.
\qed

\bigskip
\section{Algebraic nets}
\bigskip
In this section, we introduce the main notion of the paper and prove the main
results. The following definition is a discrete analogue of the definition of
algebraic webs (cf. \cite{Bl}). 
\begin{definition}
Let $N=(\A,\X)$ be a 3-net in $\P^2$ and $N^*=(\A^*,\X^*)$ be the pair of dual sets
of points and lines respectively in the dual plane $\P^*$.
We say that $N$ is {\it algebraic} if there exists a cubic curve
$\C\subset \PP^*$ containing $\A^*$ in the set $\C^0$ of its regular points.
\end{definition}

Notice that the cubic in the above definition is not assumed smooth nor even
irreducible.

In order to start proving that some nets are algebraic we need an auxiliary
definition.
\begin{definition} \label{def_complete}
We call a subset of 9 points from $\P^2$ {\it complete} it is the complete
intersection of two split cubics.
\end{definition}

The famous and needed for our purposes property of complete sets is the
Chasles theorem \cite{EGH}: if a cubic
passes through arbitrary 8 points of a complete set $T$ then it passes through
all 9 of them. Moreover the cubics passing through a complete set form a
pencil, a one-dimensional projective space spanned by the split cubics
$\C_1$ and $\C_2$ from the above definition. 

In the proofs of the next two results
we will consider certain subsets $T\subset\A^*$
with the classes $\A_i^*$ identified to a set $H$ under some bijections. 
To distinguish elements from different classes we 
write sets as $T=\{T_1|T_2|T_3\}$ where $T_i\subset H$ 
identified with a subset of $\A_i^*$
 (we will suppress other $\{\ \}$). Also if $a_i\in H$
($i=1,2,3$) we denote by $(a_1|a_2|a_3)$ the line from $\X^*$ passing through
the respective points from the classes.
\begin{theorem}
\label{thm_cyclic}
Let $(\A,\X)$ be a 3-net representing a cyclic group $\ZZ_m$.
Then this net is algebraic.
\end{theorem}
\proof
If $m=2,3$ the set $\A^*$ contains not more than 9 points that all lie on the
union of 3 lines away from their intersection points. The result follows.

In the rest of the proof we assume that $m>3$.
The set $\A^*$ is partitioned in 3 classes
$\A^*_i$ ($i=1,2,3$) of cardinality $m$. Our method of proof 
is to find a sufficient
supply of complete subsets of $\A^*$, three from each class, 
such that the split cubics from the definition \ref{def_complete} are unions
of lines from $\X^*$. 

Since the net realizes $\ZZ_m$ we can identify each class $\A_i^*$ with 
$\ZZ_m$ such that $a_1+a_2+a_3=0$ for $a_i\in\A_i^*$ ($i=1,2,3$)
if and only if the points are collinear (and lie on a line from $\X^*$). Also we
write elements of $\ZZ_m$ as integers defined modulo $m$.

To illustrate our notation and start the proof
consider the set $T_0=\{0,1,2|0,1,2|-1,-2,-3\}$.
It is complete. Indeed for the two cubics
$\C_1=(0|1|-1)\cup(1|2|-3)\cup(2|0|-2)$ and
$\C_2=(0|2|-2)\cup(1|0|-1)\cup(2|1|-3)$ we clearly have $\C_1\cap\C_2=T_0$.
We denote by $\C$ the cubic from the pencil with the base $T_0$ containing also
the point $3\in\A_1$.  It follows immediately from B\'ezout's theorem that 
$T_0\subset \C^0$.

Now we note that if $\{A_1|A_2|A_3\}$ is a complete set then $\{A_1+a_1|A_2
+a_2|A_3-a_1-a_2\}$ is also complete. This gives us several complete sets
obtained from $T_0$: $T_{k,l}=\{k,k+1.k+2|l,l+1,l+2|-k-l-1,
-k-l-2,-k-l-3\}$, where $k,l=0,1,\ldots, m-1$. 

Now we introduce sets $T(s)$ of $3(s+1)$ points, $s=2,3,\ldots,m-1$, putting
$T(s)=\{0,\ldots,s|0,\ldots,s|1,\ldots,s+1\}$.
Since $T(m-1)=\A^*$ it suffices to prove that $T(s)\subset\C^0$ for all
$s=2,\ldots,m-1$. We use finite induction on $s$. 

First $T(2)=T_0\subset\C^0$ by construction. Suppose $s>2$ and $T(s-1)\subset\C^0$.
The inductive step consists actually of three steps. First we notice that
8 points of $T_{s-3,1}$ are on $\C^0$ whence the ninth one ($s+1$ from $\A_3$)
is on $\C^0$ also. Now 8 points of $T_{s-2,0}$ are on $\C^0$ whence $\C^0$ 
also passes
through $s$ from $\A_3$. Finally the same is true for $T_{0,s-2}$ which implies
$T(s)\subset\C^0$. (The last step is redundant for $s=3$).
\qed

\medskip
Now we are ready to prove our main theorem about the algebraic character of nets
realizing finite Abelian groups.
\begin{theorem}
\label{thm_oncubic}
Let $H$ is a finite Abelian group with at least one element of order
greater than 9. Then every realization of $H$ by a 3-net in $\P^2$ is algebraic.
\end{theorem}
\proof
It follows from the condition on $H$
that $H=\ZZ_m\oplus G$ where $m>9$ and $G$ is an
Abelian group, say of order $d$.
In particular every class $\A_i$ can be partitioned into $d$ parts
of size $m$ each and this partition induces a partition of the net into
$d^2$ subnets representing each the group $\ZZ_m$. Due to Theorem \ref{thm_cyclic},
for each of those subnets there exists a cubic in the dual plane
containing all the $3m$ points dual to the lines of the net. Clearly there is a
linear ordering on these cubics such that every two subsequent cubics have
at least $m$ points in common. 

Now the proof branches out. First suppose that at least one of these cubics is
irreducible. Then since $m>9$ the cubics coincide with each other which
concludes the proof in this case. 

Second suppose that all the cubics are reducible and at least one has a line,
say $L$, and an irreducible quadric, say $Q$, as its components.  It follows from 
the collinearity pattern of the points that the only possible way in which
they can be distributed
among the components is such that the points from one class, say $\A_1^*$,
lie on $L$ and others on $Q$. Similarly to the previous 
paragraph we deduce that the quadric $Q$ is the common component of all these
cubics (using induction on the linear order and the inequality $m>6$).
In order to prove that the line components also coincide we construct more
cubics. First notice that any 9-tuple $T=\{a,b,c|a,b,c|a+b,a+c,b+c\}$ is complete
for every subset $\{a,b,c\}\subset H$. Indeed $T$ is the intersection of the
cubics $\C_1=(a|b|a+b)\cup(c|a|a+c)\cup(b|c|b+c)$
and $\C_2=(b|a|a+b)\cup(a|c|a+c)\cup(c|b|b+c)$. Notice also that any cubic
of the pencil generated by $\C_1$ and $\C_2$ has 6 points lying on $Q$ whence we
can fix a cubic $\C(T)$ from this pencil having 7 such points. Thus $\C(T)=Q\cup
L(T)$ where $L(T)$ is it linear component.

Now fix two subnets representing $\ZZ_m$ lying on cubics with linear components
$L_1$ and $L_2$ respectively. Line  $L_i$ contains a subset $\A_{1i}\subset \A_1$
of $m$ points. Choose a 9-tuple $T_1$ as in the previous paragraph with the
extra condition $a,b\in\A_{11}$ and $c\in\A_{12}$. Choose also a 9-tuple $T_2$
corresponding to a subset $\{a,b',c\}\subset H$ 
with $b'\in\A_{12}$. Then $L(T_1)$ and $L(T_2)$
have the points $a$ and $c$ from $\A_1$ in common whence coincide.
On the other hand $L(T_i)$ has 2 points from $\A_1$ common with $L_i$ whence 
$L_1=L_2$. This implies all the linear components of the cubics
containing the subnets coincide and hence the cubics themselves coincide.

Finally suppose that all the cubics containing subnets split into union of
three lines each. This means that for the dual $(3,m)$-nets each class is a
pencil. Then by Corollary \ref{cor_pencils} the whole initial net has this
property whence the dual net lies on the split cubic formed by the lines dual to
the bases of the pencils. This completes the proof.
                \qed

\medskip
\begin{remark}
(ii) Notice that 
the condition $m>6$ suffices for the conclusion of Theorem \ref{thm_oncubic}
in the case where a class of the net is a pencil. The case where all three
classes are pencils needs only $m>2$. In fact, due to Theorem \ref{Z_2},
the only unknown case in the latter
situation is for the group $\ZZ^2_2$.

(ii) Notice also that the commutativity of $H$ was really used only in the
second case where each cubic containing a subnet had two irreducible
components.
\end{remark} 

\medskip
When a $3$-net is algebraic it realizes a subgroup of the group defined on the
smooth points of the cubic by the collinearity relation.
First we recall these groups.

 If a cubic $\C$ is smooth than this 
relation on $\C$
uniquely defines the third point $\mu(a,b)\in C$ for any given $a,b\in C$.
Taking an arbitrary inflection point $p_0\in C$ as 0 one gets inverses via
$-a=\mu(0,a)$ and the binary operation defined by $a+b=-\mu(a,b)$. This
operation gives an Abelian group isomorphic to 2-torus. Three points
$p,q,r\in \C$ are collinear if and only if $p+q+r=0$.

If a cubic $\C$ is not smooth we denote as above by $\C^0$ the complement in $\C$
to the singular locus.
The following results about the algebraic structures on $C^0$ are probably
very old but I could not find them in the literature.

In the following proposition we always denote by $\CC^*$ the multiplicative group
on nonzero complex numbers and by $\CC$ the additive group on all complex numbers.
\begin{proposition}
\label{prop_groups}
(1). Let $\C$ be an irreducible singular cubic. Define the operation *
on $\C^0$ in exactly the same way as it was done on a smooth cubic. Then $\C^0$
is invariant with respect to * and forms a group isomorphic to $\CC^*$ or $\CC$
depending on the singular point being a node or a cusp respectively.

(2). Let $\C$ have as the components a quadric $Q$ and a line $L$. Then $\C^0$ is a
disjoint union of $Q^0$ and $L^0$ that are the complements in $Q$ and $L$ 
respectively to their common points. 
The collinearity relation gives
a pairing $Q^0\times Q^0\to L^0$. There exist
isomorphisms of $Q^0$ and $L^0$ (as quasiaffine varieties)
to either $\CC^*$ or $\CC$
depending on $|L\cap Q|=2$ or $1$ respectively.

(3). Finally let $\C$ be the union of three distinct
lines $L_1$. $L_2$, and $L_3$. Denote
by $L_{i}^0$ the complement in $L_i$ to the intersection points with the other
two lines, i.e., $\C^0=\cup_iL_{i}^0$. Then the collinearity relation defines a
pairing $L_1^0\times L_2^0\to L_3^0$. There exist isomorphisms
of $L_i^0$ to either $\CC^*$ or $\CC$ 
depending on $|\cap_iL_i|=3$ or 1 respectively. 
\end{proposition}
\proof
That $C^0$ is invariant with respect to the operation follows
for all the cases from B\`ezout's theorem. The other claims are proved by
choosing convenient parametrizations and, for (1), the neutral element. We will 
describe the choices leaving the straightforward checking to the reader.

(1a) Suppose that $\C$ is irreducible and has a node. We can assume that its
equation is $x^3+x^2z-y^2z=0$. A well-known parametrization of $C$ over $\P^1$
is given by $(u_1:u_2)\mapsto (u_2(u_1^2-u_2^2):u_1(u_1^2-u_2^2):u_2^3)$.
This parametrization induces an isomorphism of $\CC^*\cong\P^1\setminus\{(1:\pm
1)\}$ to $\C^0=\C\setminus \{(0:0:1)\}$. Putting $t=\frac{u_1-u_2}{u_1+u_2}$
we get this isomorphism in the form $t\mapsto \phi (t)=(4t(1-t):4t(1-t):(1-t)^3)$
($t\in\CC^*,\phi(t)\in \C^0$). 
Then $\phi$ sends the multiplication on $\CC^*$ to the
operation on $\C^0$ defined by the collinearity and the neutral element
$\phi(1)=(0:1:0)$.

(1b) Suppose that $\C$ is irreducible and has a cusp. The standard equation is
$x^3-y^2z=0$ and an isomorphism $\phi:\CC\to \C^0$ is given by 
 $\phi(s)=(s:1:s^3)$. This isomorphism sends the addition on $\CC$ to the
operation on $\C^0$ defined by the collinearity and the neutral element 
$\phi(0)=(0:1:0)$.

(2a) Suppose that the irreducible components of $\C$ are a quadric $Q$ and a line
$L$ and $|Q\cap L|=2$. The standard equations are $x^2+y^2-z^2$ and $z=0$
respectively. Isomorphisms $\phi_1:\CC^*\to Q^0$ and $\phi_2:\CC^*\to L^0$ are
given by $\phi_1(s)=(1+s^2:i(1-s^2):2s)$ and $\phi_2(s)=(1-s:i(1+s):0)$
respectively. These isomorphisms send the multiplication on $\CC^*$ to the
pairing defined by the collinearity. 

(2b) Suppose the irreducible components of $\C$ are like in the previous case
but $L$ is tangent to $Q$. The standard equations are $x^2-yz=0$ and $z=0$.
Isomorphisms $\CC\to Q^0$ and $\CC\to L^0$ are given by $s\mapsto (s:s^2:-1)$
and $s\mapsto (1:s:0)$. They send the addition to the pairing defined by the
collinearity.

(3a) Suppose that the components of $\C$ are distinct lines $L_i$ ($(i=1,2,3$)
not passing all through a common point. They can be described
by the equations $x=0$, $y=0$, and
$z=0$ respectively. Isomorphisms $\CC^*\to L_i$ are given by $s\mapsto(0:s:1)$,
$s\mapsto(0:1:s)$, and $s\mapsto (1/s:-s:0)$ respectively. They send the
multiplication to the pairing defined by the collinearity.

(3b) Finally suppose that the components of $\C$ are distinct lines $L_i$
passing through a point $p$. They can be given by the equations $x=0$, $y=0$,
and $x-y=0$ respectively. Isomorphisms $\CC\to L_i$ are given by $s\mapsto
(0:1:s)$, $s\mapsto(1:0:s)$, and $s\mapsto (1:1:s)$. They send the addition to
the needed paring.                              \qed

From this Proposition and Theorem \ref{thm_oncubic} we have the following
immediate corollary.
\begin{corollary}
\label{nonrealization}
Let $H$ be a finite Abelian group.

(i) Suppose $H$ has an element of order $\geq 10$. Then $H$ can be realized by a
3-net in $\PP^2$ if and only if $H$ has at most two invariant factors.

(ii) Suppose that $H$ has an element of order $\geq 7$. Then $H$ can be realized
by a 3-net one class of which is a pencil if and only if $H$ is cyclic.

\end{corollary}
 \bigskip
\section{Conjectures and open problems}
\bigskip

We conjecture that the conditions on the order of elements of $H$ in Corollary 
\ref{nonrealization} are not
necessary. Due to Theorem \ref{cubic}, we know that it is true in one direction.

\begin{conjecture}
\label{conj_groups}
If a finite group $H$ can be realized by a 3-net in $\PP^2$
 then $H$ is a direct sum of at
most two cyclic groups. If one class of the net is a pencil then $H$ is cyclic.
\end{conjecture}

Using the result for $\ZZ^3_2$, to settle this conjecture it suffices 
to prove or disprove that the groups $\ZZ_p^3$ are not realizable for 
$p=3,5,7$.

The results of the previous section lead naturally to the following problems.

{\bf Problem 1}. Do there exist 4-nets in $\PP^2$ nonisomorphic to the Hessian
configuration?

\medskip
{\bf Problem 2}. Do there exist 5-nets in $\PP^2$?

\medskip
{\bf Problem 3}. Does there exist a 3-net in $\PP^2$ not realizing an Abelian
group?

\medskip
{\bf Problem 4}. Do there exist nonalgebraic 3-nets in $\PP^2$?

\medskip
These problems are related. For instance the negative solution of Problem 4
would solve Problem 2 completely (and negatively)
and make a significant progress in Problem 1.

\end{document}